\documentclass{elsart}

\usepackage{amsmath,amsfonts}

\newenvironment{proofref}[1]{\begingroup\xdef\whichprop{#1}\begin{pf}}{\end{pf}\endgroup}

\numberwithin{thm}{section}

\newcommand{\dg}{d}
\newcommand{\Circ}{\mbox{\rm{Circ}}}
\newcommand{\Z}{\mathbb{Z}}

\newcommand{\arc}[1]{\mathrel{\mathop{\hbox{\
\vrule height 3pt depth -2.25pt width 11pt \hskip 1pt \
}}\limits^{\textstyle#1}}}

\newcommand{\path}[1]{\mathrel{\mathop{\hbox{\
\vrule height 3pt depth -2.25pt width 3pt \hskip 1.5pt 
\vrule height 3pt depth -2.25pt width 3pt \hskip 1.5pt 
\vrule height 3pt depth -2.25pt width 3pt \hskip 1.5pt 
\vrule height 3pt depth -2.25pt width 3pt \hskip 1.5pt 
\vrule height 3pt depth
-2.25pt width 3pt \ }}\limits^{\textstyle#1}}}

\theoremstyle{definition}
\newtheorem{notation}[thm]{Notation}
\newtheorem{assump}[thm]{Assumption}

\newcounter{subcase}

\renewcommand{\thesubcase}{\arabic{case}.\arabic{subcase}}

\newcommand{\pref}[1]{{\rm(}\ref{#1}\/{\rm)}}

\begin{document}

\begin{frontmatter}

\title{Hamiltonian cycles in $(2,3,c)$-circulant digraphs}

\author{Dave Witte Morris},
\ead{Dave.Morris@uleth.ca}
\author{Joy Morris},
\ead{Joy.Morris@uleth.ca}
\author{Kerri Webb}
\ead{Kerri.Webb@uleth.ca}

\address{Department of Mathematics and Computer Science,
\\
University of Lethbridge,
Lethbridge, Alberta, T1K~3M4, Canada}

\begin{abstract}
Let $D$ be the circulant digraph with $n$~vertices and connection set $\{2,3,c\}$. (Assume $D$ is loopless and has outdegree~$3$.) Work of S.\,C.\,Locke and D.\,Witte implies that if $n$~is a multiple of~$6$, $c \in \{(n/2) + 2, (n/2) + 3\}$, and $c$~is even, then $D$ does \emph{not} have a hamiltonian cycle. 
For all other cases, we construct a hamiltonian cycle in $D$.
\end{abstract}

\end{frontmatter}

\section{Introduction}\label{S:intro}

For $S \subset \Z$, the \emph{circulant digraph} with vertex set \/~$\Z_n$ 
and arcs from~$v$ to $v + s$ for each $v \in \Z_n$ and $s \in S$ is
denoted $\Circ(n; S)$. A fundamental open problem is to determine
which circulant digraphs have hamiltonian cycles.
By the following elegant result, circulant digraphs of
outdegree three are the smallest digraphs that need to be considered.

\begin{thm}[{R.\,A.\,Rankin, 1948, \cite[Thm.~4]{Rankin}}] \label{T:Rankin}
The circulant digraph \\ $\Circ(n;a,b)$ 
of outdegree~$2$ has a hamiltonian cycle iff there exist $s,t \in \Z^+$, such that
\begin{itemize}
	\item $s + t = \gcd(n, a-b)$,
	and
	\item $\gcd(n, sa + tb) = 1$.
\end{itemize}
\end{thm}

S.\,C.\,Locke and D.\,Witte \cite{LockeWitte} found two infinite families of non-hamiltonian circulant digraphs of outdegree~$3$; one of the 
families includes the following examples.

\begin{thm}[{Locke-Witte, cf.\ \cite[Thm.~1.4]{LockeWitte}}]\label{T:LockeWitte}
 \ 
 \begin{enumerate}
 \item $\Circ(6m; 2, 3, 3m + 2)$ is not hamiltonian if and only if $m$ is even.
 \item $\Circ(6m; 2, 3, 3m + 3)$ is not hamiltonian if and only if $m$ is odd.
\end{enumerate}
\end{thm}

In this paper, we show that the above examples are the only loopless digraphs of the form $\Circ(n; 2,3, c)$ that have outdegree~$3$ and are not hamiltonian:

\begin{thm}\label{T:main}
Assume $c \not\equiv 0,2,3 \pmod{n}$.
The digraph\/ $\Circ(n; 2, 3, c)$ is \textbf{not} hamiltonian iff 
all of the following hold
\begin{enumerate}
\item $n$ is a multiple of\/~$6$, so we may write $n = 6m$,
\item either $c=3m+2$ or $c=3m+3$,
and
\item $c$ is even.
\end{enumerate}
\end{thm}

The direction ($\Leftarrow$) of Theorem~\ref{T:main} is a restatement of part of the Locke-Witte Theorem \pref{T:LockeWitte}, so we need only prove the opposite direction. 

\begin{ack}
The work of D.\,W.\,M.\ and J.\,M.\ was partially supported by research grants from Canada's National Science and Engineering Research Council.
\end{ack}

\section{Preliminaries}

Our goal is to establish Theorem~\ref{T:main}($\Rightarrow$). We will prove the contrapositive.

\begin{notation}
Let $v_1,v_2$ be vertices of\/ $\Circ(6m;2,3,c)$ and let $s \in \{2,3,c\}$.
\begin{itemize}
\item The arc from $v_1$ to $v_1+s$ is called an $s$-arc.
\item If $v_1 + s = v_2$, we use $v_1 \arc{s} v_2 $ to denote the $s$-arc from~$v_1$ to~$v_2$.
\item If $v_1 + k s = v_2$ for some natural number~$k$, we use $v_1 \path{s}  v_2$ to denote the path $v_1, v_1+s, v_1 + 2s, \ldots, v_2$.
\end{itemize}
\end{notation}

\begin{assump} \label{A:basic}
Throughout the paper:
\begin{enumerate}
\item We assume the situation of Theorem~\ref{T:main}, so $n,c \in \Z^+$, and $c \not\equiv 0,2,3 \pmod{n}$.
\item \label{Aa:not+-1} We may assume $c \not\equiv 1, -1 \pmod{n}$. {\rm(}Otherwise, $\Circ(n;2,3,c)$ has a hamiltonian cycle consisting entirely of 
$c$-arcs.{\rm)}
\item Since the vertices of $\Circ(n;2,3,c)$ are elements of $\mathbb{Z}_n$,
we may assume $3 < c < n$.
\item We assume $n$ is divisible by~$6$ and write $n=6m$. {\rm(}Otherwise, $\Circ(n;2,3,c)$ has either a hamiltonian
cycle consisting entirely of $2$-arcs or a hamiltonian
cycle consisting entirely of $3$-arcs.{\rm)}
\end{enumerate}
\end{assump}

\begin{notation}
Let $H$ be a subdigraph of $\Circ(6m; 2,3,c)$, and let $v$ be a vertex of~$H$.
\begin{enumerate}
\item We let $\dg^+_H(v)$ and $\dg^-_H (v)$
denote the number of arcs of~$H$ directed out of, and into, vertex
$v$, respectively. 
\item If $\dg^+_H(v)=1$, and the arc from~$v$ to~$v + a$ is in~$H$, then
we say that $v$ \emph{travels by~$a$ in~$H$}.
\end{enumerate}
\end{notation}

\begin{notation} 
Let $u$ and~$w$ be integers representing vertices of\/ $\Circ(6m;2,3,c)$.
If $u-1 \le w < u+n$, let 
	$$I(u,w) = \{u,u+1,\ldots,w\}$$
be the interval of vertices from $u$ to~$w$. {\rm(}Note that $I(u,u) =
 \{u\}$ and $I(u,u-1) = \emptyset$.{\rm)}
\end{notation}

We now treat two simple cases so that they will not need to be considered in later sections.

\begin{lem}\label{L:-2or-3}
For any~$m$, $\Circ(6m;2,3,6m-2)$ and \/ $\Circ(6m;2,3,$ 
$6m-3)$ 
have hamiltonian cycles.
\end{lem}

\begin{pf}
The following is a hamiltonian cycle in $\Circ(6m;2,3,6m-2)$, where we
use $-2$ to denote the $(6m-2)$-arc:
$$ \begin{matrix}
0 &\path{2} &4& \arc{3} &7& \path{-2} &3& \arc{3} &6\\
  &\path{2} &10& \arc{3} &13& \path{-2} &9& \arc{3} &12\\
  &&& \vdots \\
  &\path{2} &6m-8& \arc{3} &6m-5& \path{-2} &6m-9& \arc{3} &6m-6\\
  &\path{2} & 6m-2& \arc{3} & 1& \path{-2} & 6m-3& \arc{3} & 0
 . \end{matrix}$$
The following is a hamiltonian cycle in $\Circ(6m;2,3,6m-3)$, where
we use $-3$ to denote the $(6m-3)$-arc:
\begin{align*}
 0 &\path{3} 6m-6 \arc{2} 6m-4 \path{-3} 2 \arc{2} 4 
 \\& \path{3} 6m-5 \path{2} 1 \arc{-3} 6m-2 \arc{2} 0
  . \qed \end{align*}
\end{pf}

\section{Most cases of the proof}\label{S:c>3m}

In this section, we prove the following two results that cover most of the cases of Theorem~\ref{T:main}:

\begin{prop}\label{P:cnot3}
If $c \le 3m$ and $c \not\equiv 3 \pmod{6}$,
then $\Circ(6m;2,3,c)$ has a hamiltonian cycle.
\end{prop}

\begin{prop}\label{P:bigc}
If $c>3m$ and $c \notin \{3m+2,3m+3\}$, then $\Circ (6m;2,3,c)$ has a hamiltonian cycle.
\end{prop}

\begin{notation}
For convenience, let $c' = 6m - c$.
\end{notation}

Note that $\Circ(6m;2,3,c) = \Circ(6m; 2,3, -c')$, and thus 
by Assumption~\ref{A:basic}(\ref{Aa:not+-1}) and Lemma~\ref{L:-2or-3},
we may assume $3 < c' < 6m-3$.

\begin{rem}
The use of~$c'$ is very convenient when $c$~is large {\rm(}so one should think of~$c'$ as being small --- less than $3m$\/{\rm)}, but it can also be helpful in some other cases.
\end{rem}

\begin{defn}
A subdigraph~$P$ of\/ $\Circ(6m;2,3,-c')$ is a \emph{pseudopath}
from~$u$ to~$w$ if $P$ is the disjoint union of a path from~$u$
to~$w$ and some number\/ {\rm(}perhaps~$0${\rm)} of cycles. In other words, if $v$ is a vertex of~$P$, then
\begin{equation*}
\dg^+_P(v) = \begin{cases} 0 &\text{if $v=w$;}\\
1 &\text{otherwise};\end{cases} \qquad \text{and} \qquad
\dg^-_P(v) = \begin{cases} 0 &\text{if $v=u$;}\\
1 &\text{otherwise}.\end{cases}
\end{equation*}
\end{defn}

\begin{defn} 
Let $u, w$ be integers representing vertices of $\Circ(6m;2,3,c)$.
If $u + c' + 2 \le w \le
u + 2c'$, let $P(u,w)$ be the pseudopath from~$u+1$ to $w-1$ whose
vertex set is~$I (u,w)$, such that $v$ travels by 
$$ \begin{cases}
\hfil 2, & \mbox{if $v \in I(u,w-c'-3) \cup I(u+c'+1,w-2)$}, \\
\hfil 3, & \mbox{if $v \in I(w-c'-2,u+c'-1)$,} \\
\hfil -c', & \mbox{if $v \in \{u + c', w\}$.}
\end{cases}$$
Notice that the range of values for $w$  because $c'>3$.
\end{defn}

\begin{lem}\label{L:P(a,b)path}
$P(u,w)$ is a  path if any of the following hold:
\begin{itemize}
\item $w - u \equiv 2c' \pmod{3}$; or 
\item $w - u \equiv 2c' +1
\pmod{3}$ and $w - u \equiv c' \pmod{2}$; or 
\item $w - u \equiv 2c'
+2 \pmod{3}$ and $w - u \not\equiv c' \pmod{2}$.
\end{itemize}
\end{lem}

\begin{pf}
We may assume that $u = 0$. 
Let $\varepsilon \in \{1,2\}$ be
 such that $w - c' -  \varepsilon
- 1 $ is even.

When $w  \equiv 2c' \pmod{3}$,
the path in $P(0,w)$ is
\begin{align*}
1 &\path{2}  w - c' - \varepsilon \path{3} c' - \varepsilon + 3 \path
{2} w \arc{-c'} w - c' 
\\& \path{3} c' \arc{-c'} 0 \path{2} w - c' +
\varepsilon - 3 \path{3} c' + \varepsilon \path{2} w - 1. 
\end{align*}

When $w \not \equiv 2c' \pmod{3}$, the hypothesis of the Lemma
implies that $w \equiv 2c' + \varepsilon \pmod{3}$. In this case, 
the path in $P(0,w)$ is
\begin{align*}
1 &\path{2}  w - c' - \varepsilon \path{3} c' \arc{-c'} 0 \path{2} w -
c' + \varepsilon - 3 
\\& \path{3} c' - \varepsilon + 3 \path{2} w
\arc{-c'} w - c' \path{3} c' + \varepsilon \path{2} w - 1. \end{align*} 
In both cases, it can be verified that the path from $1$ to $w-1$ 
contains all vertices in $I(0,w)$, and thus $P(0,w)$ is a path. (Note that
it suffices to check that the path contains both $c$-arcs, for then
$P(0,w)$ cannot contain any cycles.)
\qed \end{pf}




\begin{lem} \label{I(a,b)HamPath}
Let $k \in \Z$ be such that
\begin{itemize}
\item $k \le 6m$,
\item $c'+3 \le k \le 2c'+2$, and
 \item $k + c' \not\equiv 3 \pmod{6}$.
\end{itemize}
Let $u, w$ be integers representing vertices of $\Circ(6m;2,3,c)$.
If $u \leq w$ and $w-u+1=k$, then the subgraph
induced by $I(u,w)$ has a hamiltonian path that starts at $u+1$
and ends in $\{w-1,w\}$.
\end{lem}

\begin{pf}
We consider three cases.

\setcounter{case}{0}

\begin{case}
Assume $k \equiv 2c' + 1 \pmod{3}$.
\end{case}

We have $w - u = k - 1 \equiv 2c' \pmod{3}$. Since $k  \leq 2c'+2$,
this implies $w-u \leq 2c'$.
By
Lemma~\ref{L:P(a,b)path}, $P(u,w)$ is a hamiltonian path from
$u+1$ to $w-1$.

\begin{case}
Assume $k \equiv 2c' + 2 \pmod{3}$.
\end{case}

Suppose, first, that $k \neq c' + 3$ (so $k \ge c' + 4$). Letting
$w' = w - 1$, then
$$ w' - u = w - u - 1 = k - 2 \ge (c'+4) - 2 = c' + 2 $$
and $w' - u = (w - 1) - u = k - 2 \equiv 2c' \pmod{3}$. By
Lemma~\ref{L:P(a,b)path}, $P(u,w')$ is a hamiltonian path in
$I(u,w')$ from $u+1$ to $w' - 1$. Adding the $2$-arc from $w' - 1$
to $w' + 1 = w$ yields a hamiltonian path in $I (u,w)$ from $u+1$
to~$w$.

Suppose instead that $k = c' + 3$. Then $w - u = k - 1 \equiv 2c' +
1 \pmod{3}$ and $w - u = k - 1 = (c' + 3) - 1 \equiv c' \pmod{2}$,
so by Lemma~\ref{L:P(a,b)path}, $P(u,w)$ is a hamiltonian path from $u+1$ to $w-1$.

\begin{case}
Assume $k \equiv 2c' \pmod{3}$.
\end{case}

By assumption, we have $k + c' \equiv 2c'+c' \equiv 0\pmod{3}$.  Since $k+c' \not\equiv 3 \pmod{6}$,
we must have $k+c' \equiv 0\pmod{6}$, so $k \equiv c' \pmod{2}$. Then $w - u = k - 1
\equiv 2c' + 2 \pmod{3}$
   and $w - u \equiv k - 1 \not\equiv k \equiv c' \pmod{2}$,
so by Lemma~\ref{L:P(a,b)path}, $P(u,w)$ is a hamiltonian path from $u+1$ to $w-1$.
\qed \end{pf}

It is now easy to prove Propositions~\ref{P:cnot3} and~\ref{P:bigc}.

\begin{proofref}{P:cnot3}
As previously mentioned, we may assume $c>3$.
Since $3 < c \le 3m$, we have $3m \le c' < 6m - 3$, so 
	$$c' + 3 < 6m < 2c' + 2 .$$
Furthermore, since $c \not\equiv 3 \pmod{6}$, we have 
	$$6m + c' \equiv c' \not\equiv 3 \pmod{6} .$$
Hence, Lemma~\ref{I(a,b)HamPath} implies that the interval $I(0, 6m-1)$ has a hamiltonian path from $1$ to $6m-2$ or to $6m-1$. Inserting the $3$-arc from $6m-2$ to~$1$ or the $2$-arc from $6m-1$ 
to~$1$, yields a hamiltonian cycle. Since $I(0,6m-1)$ is the entire digraph, this completes the proof.
 \qed \end{proofref}

 \begin{lem} \label{L:sumK}
Let $\mathcal{K}$ be the set of integers~$k$ that satisfy the
conditions of Lemma~\ref{I(a,b)HamPath}.
If\/ $4 \le c' < 3m$, and $n_0 \ge 2(c' + 4)$, then either
	\begin{enumerate}
	\item $n_0$ can be written as a sum $n_0 = k_1 +
k_2 + \cdots + k_s$, with each $k_i \in \mathcal{K}$,
	or
	\item $c' = 6$ and $n_0 = 29$.
	\end{enumerate}
\end{lem}

\begin{pf}
Note that, since $c' < 3m$, we have $2c'+2 \le 6m$, so the first inequality in the definition of~$\mathcal{K}$ is redundant --- it can be ignored.

Let us treat some small cases individually:
\begin{itemize}
\item If $c' = 4$, then $\mathcal{K} = \{7,8,9,10\}$. It is easy to see that every integer $\ge 14$ is a sum of elements of~$\mathcal{K}$.
\item If $c' = 5$, then $\mathcal{K} = \{8,9,11,12\}$. It is easy to see that every integer $\ge 16$ is a sum of elements of~$\mathcal{K}$.
\item If $c' = 6$, then $\mathcal{K} = \{10,11,12,13,14\}$. It is easy to see that every integer $\ge 20$ is a sum of elements of~$\mathcal{K}$, except that $29$ is not such a sum.
\end{itemize}

Henceforth, we assume $c' \ge 7$, so $4c' + 5 \ge 3(c'+4)$. Then, since $c' + 4 \in \mathcal{K}$, we may assume, by subtracting some multiple of $c' + 4$, that 
	$$n_0 \le 3(c'+4) - 1 \le 4c' + 4 .$$
Under this assumption, we prove the more precise statement that 
$$ \text{$n_0 = k_1 + k_2$, with $k_1, k_2 \in \mathcal{K}$, and $k_1 \le k_2 \le k_1 + 3$.} $$
Assume that $n_0$ cannot be written as such a sum. (This will lead to a contradiction.)
Because $2(c'+4) = (c'+ 4) + (c'+4)$, we must have $n_0 > 2(c'+4)$. 
Then, by induction, we may assume there exist $k_1,k_2 \in \mathcal{K}$, such that $k_1 \le k_2 \le k_1 + 3$ and
	$$k_1 + k_2 = n_0 - 1 .$$
Since $2c' + 1, 2c' + 2 \in \mathcal{K}$, we must have 
	$$n_0 \le 4c' + 1 ,$$
so 
	$$k_1 \le 4c'/2 = 2c' .$$
Note that
	\begin{align*}
	n_0	&= (k_1 + 1) + k_2, \\
	n_0	&= (k_1 + 2) + (k_2 - 1), \\
	n_0	&= k_1 + (k_2 +1) .
	\end{align*}
From the first equation (and the fact that $k_2 \in \mathcal{K}$, we see that $k_1 + 1 \notin \mathcal{K}$. 
Because $\mathcal{K}$ contains $5$ of any~$6$ consecutive integers between 
$c'+3$ and $2c'+2$, this implies that $k_1 + 2 \in \mathcal{K}$. Hence, the second equation implies $k_2 - 1 \notin \mathcal{K}$. Therefore $k_2 - 1 = k_1 + 1$, so
	$$ k_2 + 1 = k_1 + 3 .$$
Hence, the third equation implies $k_2 + 1 \notin \mathcal{K}$, so we must have $k_2 + 1 > 2c' + 2$; therefore $k_2 = 2c' + 2$, which implies $k_1 = 2c'$. Hence
	$$ 4c' + 2 = k_1 + k_2 = n_0 - 1 \le 4c' .$$
This is a contradiction.
\qed \end{pf}

\begin{proofref}{P:bigc}
Note that  $c' < 3m$.
We dealt with the cases $c=6m-2$ and $c=6m-3$
in Lemma~\ref{L:-2or-3}, and the cases $6m \in \{2c'+4, 2c'+6\}$ are dealt with by Theorem~\ref{T:LockeWitte}.  Furthermore, we noted in Assumption~\ref{A:basic} that the case $c = 6m-1$ is clearly hamiltonian. Therefore,
we may assume in what follows that $c' \ge 4$ and $6m \notin \{ 2c'+4, 2c'+6\}$.

Let $\mathcal{K}$ be the set of integers~$k$ that satisfy the
conditions of Lemma~\ref{I(a,b)HamPath}.
We claim that $6m$ can be written as a sum $6m = k_1 +
k_2 + \cdots + k_s$, with each $k_i \in \mathcal{K}$. 
If $6m \ge 2(c'+4)$, then this is immediate from Lemma~\ref{L:sumK} 
(and the fact that $6m \neq 29$). On the other hand, 
if $6m < 2(c'+4)$, then, since $6m$~is even, and $6m \notin \{2c'+4, 2c'+6\}$, 
we see that $6m=2c'+2 \in
\mathcal{K}$, so $6m$ is obviously a sum of elements
of~$\mathcal{K}$. This completes the proof of the claim. 

The preceding paragraph implies that we may
cover the vertices of $\Circ(6m;2,3,c)$ by a disjoint collection
of intervals $I(u_i,w_i)$, such that the number of vertices in
$I(u_i,w_i)$ is~$k_i$. By listing the intervals in their natural
order, we may assume $u_{i +1} = w_{i} + 1$. By
Proposition~\ref{I(a,b)HamPath}, the vertices of $I(u_i,w_i)$ can
be covered by a path~$P_i$ that starts at $u_i + 1$ and ends in
$\{w_i - 1, w_i\}$. Since 
	$$(u_{i+1} + 1) - w_i = (w_i + 2) - w_i = 2$$
and 
	$$(u_{i+1} + 1) - (w_i - 1) = (w_i + 2) - (w_i - 1) = 3 ,$$
there is an arc from the terminal vertex of~$P_i$ to the initial
vertex of~$P_{i+1}$. Thus, by adding a number of $2$-arcs and/or
$3$-arcs, we may join all of the paths $P_1,P_2,\ldots, P_s$ into
a single cycle that covers all of the vertices of
$\Circ(6m;2,3,c)$. Thus, we have constructed a hamiltonian cycle.
\qed \end{proofref}

\section{The remaining cases}\label{S:cequiv3}

In this section, we prove the following result. 

\begin{prop} \label{P:c=3mod6}
If\/ $3 < c \le 3m$ and $c\equiv 3 \pmod 6$, then $\Circ(6m; 2,3,c)$ has a hamiltonian cycle. 
\end{prop}

Proposition~\ref{P:c=3mod6} together with 
Propositions~\ref{P:bigc} and~\ref{P:cnot3} (and Theorem~\ref{T:LockeWitte}) completes the proof of Theorem~\ref{T:main}.

\begin{defn}
Let $t$ be any natural number, such that $0 \leq 6t \leq c-9$.
\begin{enumerate}
\item Let 
	\begin{align*}
	\ell_1 &= c-5, \\
	\ell_2 &= \ell_2(t) = c - 1 + 6t , \\
	\ell_3 &= c - 2, \\
	\ell_4 &= c+3 .
	\end{align*}
\item Define subdigraphs $Q_1, Q_2, Q_3$ and~$Q_4$ of $\Circ(6m; 2,3 ,c)$ as follows:
\begin{itemize}
\item The vertex set of $Q_i$ is $I(0,\ell_i + 2) \cup \{\ell_i + 5\}$.
\item In~$Q_1$, vertex~$v$ travels by
$\begin{cases}    c, &\text{if $v=0$;}\\
                                2, &\text{if $v=1$ or $2$;}\\
                                3, &\text{if $v=3,4,\ldots,c-6$.}
                \end{cases}$
\item In~$Q_2$, vertex~$v$ travels by
$\begin{cases}    c, &\text{if $v=1$ or $6t+4$;}\\
                                2, &\text{if $v=2$ or $6t+5 \leq v\leq c-2$;}\\
                                3, &\text{if $v=0$ or $3 \leq v \leq 6t+3$ or $c-1 \leq v \leq c-2 + 6t$.}
                \end{cases}$
\item In~$Q_3$, vertex~$v$ travels by
$\begin{cases}    c, &\text{if $v=1$ or $3$;}\\
                                2, &\text{if $v=1,2$ or $4 \leq v \leq c-3$.}
                \end{cases}$
\item In~$Q_4$, vertex~$v$ travels by
$\begin{cases}    c, &\text{if $v=2$ or $8$;}\\
                                2, &\text{if $9 \leq v \leq c-1$;}\\
                                3, &\text{if $v=0,1$ or $3 \leq v
                                \leq 7$ or $c \leq v \leq c+2$.}
                \end{cases}$
\end{itemize}
\end{enumerate}
\end{defn}

\begin{notation}
For ease of later referral, we also let $\ell_i(t)$ denote $\ell_i$ for $i \in \{1,3,4\}$.
\end{notation}

\begin{lem}$\;$
\begin{enumerate}
\item Each $Q_i$ is the union of four disjoint paths 
from $\{0,1,2,5\}$ to $\{\ell_i, \ell_i+1,\ell_i+2,\ell_i+5\}$.

\item More precisely, let $u_1 = 0$, $u_2 = 1$, $u_3 = 2$, and $u_4 = 5$, and 
define permutations
$$ \mbox{$\sigma_1=(1423)$, $\sigma_2=(234)$,
$\sigma_3=(1324)$, and $\sigma_4 = \mbox{\rm identity}.$}$$
Then $Q_i$ contains a path
from $u_k$ to $\ell_i + u_{\sigma_i(k)}$, for $k = 1,2,3,4$.
\end{enumerate}
\end{lem}

\begin{pf}
The paths in~$Q_1$ are:
	\begin{align*}
	&0 \arc{c} c \quad (= \ell_1 + 5) ,\\
	&1 \arc{2} 3 \path{3} c - 3 \quad (= \ell_1 + 2),\\
	&2 \arc{2} 4 \path{3} c-5 \quad (= \ell_1),\\
	&5 \path{3} c-4 \quad (= \ell_1 + 1).
	\end{align*}
The paths in~$Q_2$ are:
	\begin{align*}
	&0 \path{3} 6t + 6 \path{2} c-1 \path{3} c-1 + 6t \quad (= \ell_2) ,\\
	&1 \arc{c} c+1 \path{3} c + 1 + 6t \quad (= \ell_2 + 2),\\
	&2 \arc{2} 4 \path{3} 6t+4  \arc{c} c + 4 + 6t \quad (= \ell_2 + 5),\\
	&5 \path{3} 6t + 5 \path{2} c \path{3} c + 6t \quad (= \ell_2 + 1).
	\end{align*}
The paths in~$Q_3$ are:
	\begin{align*}
	&0 \arc{c} c  \quad (= \ell_3 + 2) ,\\
	&1 \arc{2} 3 \arc{c} c + 3 \quad (= \ell_3 + 5),\\
	&2 \path{2} c-1 \quad (= \ell_3 + 1),\\
	&5 \path{2} c-2 \quad (= \ell_3).
	\end{align*}
The paths in~$Q_4$ are:
	\begin{align*}
	&0 \path{3} 9 \path{2} c \arc{3} c+3 \quad (= \ell_4) ,\\
	&1 \path{3} 10 \path{2} c + 1 \arc{3} c + 4 \quad (= \ell_4 + 1),\\
	&2 \arc{c} c + 2 \arc{3} c+5 \quad (= \ell_4 + 2),\\
	&5 \arc{3} 8 \arc{c} c + 8 \quad (= \ell_4 + 5)
	. \qed \end{align*}
\end{pf}

From the above lemma, we see that translating $Q_j$ by $\ell_i$ yields a disjoint union of $4$~paths whose initial vertices are precisely the terminal vertices of the paths in~$Q_i$. Hence, composing $Q_i$ with this translate of~$Q_j$ results in a disjoint union of $4$~paths: namely, a path from $u_k$ to $\ell_i + \ell_j + u_{\sigma_j \sigma_i(k)}$, for $k = 1,2,3,4$.
Continuing this reasoning leads to the following conclusion:

\begin{lem}\label{L:product}
If, for some natural number~$s$, there exist sequences 
\begin{itemize}
\item $I=(i_1,i_2, \ldots, i_s)$ with each $i_j \in \{1,2,3,4\}$, 
and 
\item $T = (t_1,t_2, \ldots, t_s)$ with $0 \leq 6t_j \leq c-9$, for each~$j$, 
\end{itemize}
such that
\begin{enumerate} \renewcommand{\theenumi}{\roman{enumi}}
\item $\sum_{j=1}^s \ell_{i_j}(t_j) = 6m$, and
\item the permutation product $\sigma_{i_s}\sigma_{i_{s-1}}
\cdots \sigma_{i_1}$ is a cycle of length $4$,
\end{enumerate}
then $\Circ(6m; 2,3,c)$ has a hamiltonian cycle constructed by 
concatenating $s$~appropriate translates of $Q_1$, $Q_2$, $Q_3$, and/or~$Q_4$.
\end{lem}

\begin{proofref}{P:c=3mod6}
Since $\sigma_4$ is the identity and $\ell_4=c+3$, 
we see that if $\Circ(6m;2,3,c)$ has a hamiltonian cycle
constructed by concatenating translates of $Q_1,Q_2,Q_3$, and $Q_4$,
then $\Circ(6m+c+3; 2,3,c)$ also has such a
hamiltonian cycle. Thus, by subtracting some multiple of $c+3$ from~$6m$, we may assume 
	$$2c-6 \leq 6m \leq 3c-9 .$$
(For this modified $c$, it is possible that $c>3m$.)

Recall that $0 \le 6t \le c-9$, so $2c - 6 + 6t$ can be any multiple of~$6$ between $2c - 6$ and $3c - 15$. Since
$\sigma_2 \sigma_1=(1243)$ and $\ell_1 + \ell_2(t) = 2c - 6 + 6t$, it follows that $\Circ(6m;2,3,c)$
has a hamiltonian cycle constructed by concatenating~$Q_1$ with a translate of~$Q_2$ whenever
$2c-6 \leq 6m \leq 3c-15$.

The only case that remains is when $6m = 3c - 9$. Now
$\sigma_3^2\sigma_1=(1324)$ and $\ell_1+2\ell_3 =3c-9$, so $\Circ(3c-9;2,3,c)$
has a hamiltonian cycle constructed by concatenating~$Q_1$ with two translates of~$Q_3$. 
\qed \end{proofref}

\end{document}